\documentclass{amsart}

\usepackage{amssymb,amsmath}

\numberwithin{equation}{section}

\newtheorem{theorem}{Theorem}

\newtheorem{corollary}[theorem]{Corollary}
\newtheorem{prop}{Proposition}

\newcommand\pd{{\partial}}

\newcommand{\R}{{\mathbb{R}}}
\newcommand{\C}{{\mathbb{C}}}
\newcommand{\N}{{\mathbb{N}}}

\newcommand{\h}{\mathcal{H}}
\newcommand{\cl}{\mathcal{L}}

\newcommand\cT{{\mathcal T}}
\newcommand\T{{\mathbb T}}

\newcommand\ep{\epsilon}

\newcommand\yinf{y\rightarrow+\infty}

\newcommand\limy{\lim\limits_{y\rightarrow+\infty}}
\newcommand\limz{\lim\limits_{y\rightarrow+\infty}}

\begin{document}

\title[]
{Ill-posedness of the Prandtl equations in Sobolev spaces around a shear flow with general decay}

\author{Cheng-Jie Liu}
\address{Cheng-Jie Liu
\newline\indent
Department of Mathematics, The Chinese University of Hong Kong,
\newline\indent Hong Kong, P. R.  China}
\email{cjliusjtu@gmail.com}

\author{Tong Yang}
\address{Tong Yang
\newline\indent
Department of Mathematics, Shanghai Jiao Tong University,
\newline\indent Shanghai, P. R.  China
\newline\indent
and Department of mathematics, City University of Hong Kong
\newline\indent
Hong Kong, P. R. China}
\email{matyang@cityu.edu.hk}


\subjclass[2000]{35M13, 35Q35, 76D10, 76D03, 76N20}

\date{}

\keywords{Prandtl equations, ill-posedness, linear instability, shear flow, 
monotonic velocity fields.}

\begin{abstract}

Motivated by the paper \cite{GV-D} [JAMS, 2010] about
the linear ill-posedness  for the Prandtl equations around a shear flow with exponential decay in normal variable, and the recent
study of well-posedness on the Prandtl equations in Sobolev spaces,
this paper aims
to  extend the result in \cite{GV-D} to the case when the shear flow 
has general  decay. The key observation is to
construct  an  approximate solution that captures
the initial layer to the linearized problem motivated by the 
precise formulation of solutions to the inviscid Prandtl
equations. 

\end{abstract}

\maketitle



\section{Introduction and main results}

The Prandtl equations were introduced by Ludwing Prandtl \cite{Pra} in 1904 to describe the motion of fluid  with small viscosity near a solid boundary
with non-slip boundary condition. This seminal work sets the foundation of
boundary layer theories. Even though the Prandtl equations have been proved
its importance in physics and engineering applications, the mathematical theories established are far from being satisfactory.

One of the pioneering works by Oleinik and her collaborators \cite{Ole} in 1960s
shows that under the monotonicity condition of the tangential velocity
component in the normal direction to the boundary, local well-posedness
theories of Prandtl equations can be established. This result was recently
further improved in the framework of Sobolev spaces, cf. \cite{AWXY, M-W}. 
On the other hand, the ill-posedness of this system in the Sobolev
spaces for perturbation of a shear flow with a non-degenerate critical
point was proved in the interesting paper  \cite{GV-D} linearly and
then nonlinear in \cite{GV-N,guo}, following the long time study on
the instability by many authors, cf. \cite{e-2,  grenier, LWY2, van} ect. 
It is noted that in the work  \cite{GV-D}, the shear flow is assumed
to be exponentially decay to the uniform Euler flow in the normal
direction with respect to the boundary. However, as pointed out in \cite{vicol}, the
exponential decay should not be essential, in particular, in 
the physical consideration. Therefore, it remains the question whether
the instability showed in \cite{GV-D} for exponential decay
shear flow holds  with general decay. In fact, the
answer to this question in some sense reveals the monotonicity condition
on the tangential velocity component is a necessary and sufficient
condition for well-posedness in the framework of Sobolev spaces.

In the following, we will first present the result for the
Prandtl equations in a two dimensional domain
$\Omega\triangleq\{(t,x,y): t>0,(x,y)\in\T\times\R^+\}$, and
then in the last section, we will give some discussion on the
case in three space dimensions. That is, consider
\begin{equation}\label{prandtl}
\begin{cases}
\partial_t u+u\partial_x u+v\partial_yu+\partial_xP-\partial_y^2u=0,\qquad &\\
\partial_xu+\partial_yv=0,\qquad &{\mbox in}\quad\Omega,\\
(u, v)|_{y=0}=0, \qquad \lim\limits_{y\to+\infty}u=U(t,x),
\end{cases}
\end{equation}
where  $U=U(t,x)$ and $P=P(t,x)$ are the tangential velocity and pressure 
of the Euler flow adjacent to the boundary layer. Moreover, $U(t,x)$ and $P(t,x)$ satisfy the Bernoulli equation:
$$
\partial_t U+U\pd_xU +\pd_xP=0.
$$

Since we are interested in the instability structure of this system
around a shear flow, as in \cite{GV-D}, we
consider the simple case of \eqref{prandtl} when the Euler flow $U$ is constant:
\[
U(t,x)~\equiv~U_0,\quad\mbox{and then},\quad \pd_x P(t,x)~\equiv~0.
\]
In this case, the problem \eqref{prandtl} becomes 
\begin{equation}\label{pr}
\begin{cases}
\partial_t u+u\partial_x u+v\partial_yu-\partial_y^2u=0,\qquad &\\
\partial_xu+\partial_yv=0,\qquad &{\mbox in}\quad\Omega,\\
(u, v)|_{y=0}=0, \qquad \lim\limits_{y\to+\infty}u=U_0.
\end{cases}
\end{equation}
Note that \eqref{pr} has a special shear flow solution $\big(u_s(t,y),0\big)$, where the function $u_s(t,y)$ is a smooth solution to the following heat equation:
\begin{equation}\label{heat}\begin{cases}
\pd_t u_s-\pd_y^2u_s=0,\qquad {\mbox in}\quad\Omega,\\
u_s|_{y=0}=0,\quad\limy u_s=U_0,\\
u_s|_{t=0}=U_s(y)
\end{cases}\end{equation}
with an initial shear layer $U_s(y)$.  Then, we consider the linearization of
the problem \eqref{pr} around the shear flow $\big(u_s(t,y),0\big)$, and obtain
\begin{equation}\label{linear_pr}\begin{cases}
\pd_t u+u_s\pd_xu+v\pd_yu_s-\pd_y^2u=0,\qquad &\\
\pd_xu+\pd_yv=0,\qquad &{\mbox in}\quad\Omega,\\
(u,v)|_{y=0}=0,\quad \limy u=0.
\end{cases}\end{equation}

In \cite{GV-D}, 
the authors showed that if the intial data $U_s(y)$ of
the shear flow has a non-degenerate
critical point, then the linear problem \eqref{linear_pr} is
ill-posed in the case that $u_s-U_0$  exponentially decays to zero as $\yinf$. 
The goal of this paper is to show that the exponential
decay condition is not necessary. Indeed, a physcial quantity that measures the effect of
the boundary layer matching the outer flow, called displacement thickness, cf. \cite[p.311]{batchelor}, is defined by
\begin{equation}
\delta(t,x)=\int_0^\infty (1-\frac{u(t,x,y)}{U(t,x)}) dy.
\end{equation}
Hence, the finiteness of the displacement thinkness only requires
the integrability of the above function, which admits general decay of
$u(t,x,y)$ to $U(t,x)$ when $y$ tends to infinity.

To continue, let us first introduce some notations.
Denote by $T(t,s)$ the linear solution operator:
\begin{equation}\label{def_T}
T(t,s)u_0~:=~u(t,\cdot),
\end{equation}
where $u$ is the solution to the problem \eqref{linear_pr} with $u|_{t=s}=u_0.$ Also,  for any $\alpha,m\geq0$, denote 
\[\begin{split}
&W_\alpha^{m,\infty}(\R^+)~:=~\{f=f(y),y\in\R^+;~\|f\|_{W_\alpha^{m,\infty}}\triangleq\|e^{\alpha y}f(y)\|_{W^{m,\infty}(\R^+)}<\infty\},\\
&\h_\alpha^m~:=~\{f=f(x,y),(x,y)\in\T\times\R^+;~\|f\|_{H_\alpha^m}\triangleq\|f(\cdot)\|_{H^m(\T_x,W_\alpha^{0,\infty}(\R^+_y))}<\infty\}.
\end{split}\]
The main result on the linear ill-posedness of the Prandtl equations 
can be stated as follows.

\begin{theorem}\label{thm_lin}
Let $u_s(t,y)$ be the solution of the problems \eqref{heat} satisfying
 $$u_s-U_0\in C^0\big(\R^+; W_0^{4,\infty}(\R^+)\big)\cap C^1\big(\R^+; W_0^{2,\infty}(\R^+)\big),$$
and assume that the initial shear layer $U_s(y)$ has a non-degenerate critical point in $\R^+$.
Then, there exists $\sigma>0$ such that for all $\delta>0$,
\begin{equation}\label{est_in}
\sup\limits_{0\leq s< t\leq\delta}\big\|e^{-\sigma(t-s)\sqrt{|\pd_x|}}T(t,s)\big\|_{\cl(\h_\alpha^m,\h_0^{m-\mu})}
~=~+\infty,\quad \forall \alpha,m\geq0,~\mu\in[0,\frac{1}{2}).
\end{equation}
\end{theorem}

One consequence of the above theorem gives

\begin{corollary}
Under the assumptions of Theorem \ref{thm_lin}, it holds that for any $\delta>0$ and $\alpha,m\geq0$,
 \begin{equation}\label{est_ins}
\sup\limits_{0\leq s< t\leq\delta}\big\|T(t,s)\big\|_{\cl(\h_\alpha^m,\h_0^{0})}
~=~+\infty.
\end{equation}
\end{corollary}

At the end of the introduction, let us mention that most of the mathematical
theories for the Prandtl equations before 2000 can be found in 
the excellent review article \cite{e-1}. In addition to those
works mentioned before, some other interesting works can be found
in \cite{LWY, LWY1} for the three space dimensional Prandtl
equations with special structure to avoid the secondary flow, cf. \cite{moore}, the 
works in the framework of analytic function space in \cite{cannone, S-C, Z-Z},
and 
the existence of global weak solutions in \cite{LWY1, X-Z}.

The rest of the paper will be arranged as follows. 
The main result on the  linear instability for the system around a shear
flow with general decay and a non-degenerate critical point will be
proved in the next section by a new construction of approximate solutions.
Some discussions on the case in three space dimensions
will be given in the last section.

\section{Linear instability}

In the following three subsections, we will prove Theorem \ref{thm_lin} for the linear instability of the Prandtl equations.

\subsection{Instability mechanism}

In this subsection, we firstly recall the linear ill-posedness result in \cite{GV-D} about
the linear instability mechanism of Prandtl equations, and then
introduce the new approximate solutions for general decay shear flow. The key 
observation 
in \cite{GV-D} is to construct an unstable approximate solution to \eqref{linear_pr}, in high frequency in the tangential variable $x$, with exponential growth in time $t$.
To illustrate this kind of instability mechanism,
as in \cite{GV-D}, one can  first replace the background shear flow in \eqref{linear_pr} by its initial data, and consider the following simpler problem
with coefficients independent of the variable $t$:
\begin{equation}\label{sim_lin}\begin{cases}
\pd_t u+U_s\pd_xu+vU_s'-\pd_y^2u=0,\quad&\\
\pd_xu+\pd_yv=0,\quad&{\rm in}\quad \Omega,\\
(u,v)|_{y=0}=0,\qquad\limy u=0.
\end{cases}\end{equation}
Denote by $\cl_s$ the linearized Prandtl operator in \eqref{sim_lin} around the shear flow $\big(U_s(y),0\big):$
\begin{equation}\label{def_ml}
\cl_su~:=~U_s\pd_xu+vU_s'-\pd_y^2u,\quad\mbox{with}\quad v(t,x,y)=-\int_0^y\pd_xu(t,x,z)dz.
\end{equation}

In Section 2 of \cite{GV-D}, the authors construct an approximate solution of \eqref{sim_lin}, which has high $x-$frequency
of the order  $\ep^{-1}$ and grows in $t$ exponentially at the rate of $\ep^{-\frac{1}{2}}$ for $\ep\ll1$. Precisely, one can look for solutions to \eqref{sim_lin} in the form
\[
(u,v)(t,x,y)~=~e^{i\ep^{-1}(x+w_\ep t)}\big(u_\ep(y),\ep^{-1}v_\ep(y)\big).
\]
By plugging this into \eqref{sim_lin}, the divergence free condition gives  $u_\ep(y)=iv_\ep'(y)$, and then  the first equation of \eqref{sim_lin} yields
\begin{equation}\label{pr_v}\begin{cases}
\big(w_\ep+U_s(y)\big)v_\ep'(y)-U_s'(y)v_\ep(y)+i\ep v_{\ep}^{(3)}(y)=0,\quad y>0,\\
v_\ep|_{y=0}=v_\ep'|_{y=0}=0.
\end{cases}\end{equation}
Let $a>0$ be a non-degenerate critical point of the initial shear layer $U_s(y)$, the following result was proved in \cite{GV-D}.

\begin{prop}\label{prop_GD}
 There exists an approximate solution $(u_\ep^{app},v_\ep^{app})(t,x,y)$ to the problem \eqref{sim_lin} in the form of
\begin{equation}\label{app_1}
(u_\ep^{app},v_\ep^{app})(t,x,y)~=~e^{i\ep^{-1}(x+w_\ep t)}\big(iv_\ep'(y),\ep^{-1}v_\ep(y)\big),
\end{equation}
where 
\begin{equation}\label{def_w}
w_\ep=-U_s(a)+\ep^{\frac{1}{2}}\tau
\end{equation} 
for some constant $\tau\in\C$ with the imaginary part $\Im\tau<0$, and $v_\ep(y)\in W_0^{3,\infty}(\R^+)$, such that the error term $r^{app}_\ep:=\pd_t u^{app}_\ep+\cl_s u_\ep^{app}$ satisfies
\begin{equation}\label{est_err1}
r^{app}_\ep(t,x,y)~=~e^{i\ep^{-1}(x+w_\ep t)}R^{app}_\ep(y),\quad R^{app}_\ep(y)\in W_0^{0,\infty}(\R^+).
\end{equation}
\end{prop}

In fact, as shown in \cite{GV-D},  the function $v_\ep(y)$ can be
 devided into a "regular" part $v_\ep^{reg}(y)$ and a "shear layer" part $v^{sl}_\ep(y)$, i.e.,
\begin{equation}\label{for_v}\begin{split}
v_\ep(y)&=v^{reg}_\ep(y)+v^{sl}_\ep(y)\\
&=H(y-a)\big[U_s(y)-U_s(a)+\ep^{\frac{1}{2}}\tau\big]+\ep^{\frac{1}{2}}V(\frac{y-a}{\ep^{\frac{1}{4}}}).
\end{split}\end{equation}
Here, $H$ is the Heaviside function, and
 the shear layer profile $V(z)$ solves the following ODE:
\begin{equation}\label{eq_w}\begin{cases}
\Big(\tau+U_s''(a)\frac{z^2}{2}\Big)V'-U_s''(a)zV+iV^{(3)}=0,\quad z\neq0,\\
[V]\big|_{z=0}=-\tau,\quad~[V']\big|_{z=0}=0,\quad~[V'']\big|_{z=0}=-U_s''(a),\\
\lim\limits_{z\rightarrow\pm\infty} V~=~0,\quad  exponentially,
\end{cases}\end{equation}
where the complex constant $\tau$ is the same as the one in \eqref{def_w}, and the notation $[u]\big|_{z=0}~:=~\lim\limits_{\delta_1\to 0+}u(\delta_1)-\lim\limits_{\delta_2\to 0-}u(\delta_2)$ denotes the jump of a related function $u(z)$ across $z=0$. One can check that by virtue of $ w_\ep$ given in \eqref{def_w} , the function $v_\ep(y)$ defined in
 \eqref{for_v}  solves the problem \eqref{pr_v} 
except for the $O(\ep)$-term coming from diffusion. Consequently, the corresponding approximate solution \eqref{app_1} admits the $O(\ep^{-1})$-terms of the first equation of \eqref{sim_lin}, which implies the estimate \eqref{est_err1} automatically. Indeed, the direct calculation gives the expression of the error term $R_\ep^{app}(y)$ defined in \eqref{est_err1}:
\begin{equation}\begin{split}\label{for_R}
R_\ep^{app}(y)=&-\ep^{-1}\big[U_s(y)-U_s(a)-U_s''(a)\frac{(y-a)^2}{2}\big](v_\ep^{sl})'(y)\\
&
+\ep^{-1}\big[U_s'(y)-U_s''(a)(y-a)\big]v_\ep^{sl}(y)-i(v_\ep^{reg})^{(3)}(y),
\end{split}\end{equation}
so that the estimate of $R_\ep^{app}(y)$ in \eqref{est_err1} follows from the exponential decay of the profile $V(z)$. Furthermore, we have from \eqref{for_R},
 $$R_\ep^{app}(y)+i(v_\ep^{reg})^{(3)}(y)\in W_\alpha^{0,\infty}(\R^+),~\mbox{for any }~\alpha\geq0.$$
Note that the term $i(v_\ep^{reg})^{(3)}(y)$ does not appear in the error term
when the
background profile is the shear flow $u_s(t,y)$, not the initial shear layer $U_s(y)$, because of the heat equation, cf.
\eqref{tr}.

In addition, we refer to \cite{GV-D} and note that the pair $\big(\tau,V(z)\big)$ takes the following form:\
\begin{equation}\label{ch_w}\begin{cases}
\tau~&=~\big|\frac{U_s''(a)}{2}\big|^{\frac{1}{2}}~\tilde\tau,\\
V(z)~&=~\big|\frac{U_s''(a)}{2}\big|^{\frac{1}{2}}\Big[\big(\tilde\tau+\big|\frac{U_s''(a)}{2}\big|^{\frac{1}{2}}z^2\big)
W\big(\big|\frac{U_s''(a)}{2}\big|^{\frac{1}{4}}z\big)-1_{\R^+}\big(\tilde\tau+\big|\frac{U_s''(a)}{2}\big|^{\frac{1}{2}}z^2\big)\Big],
\end{cases}\end{equation}
where the function $W(z)$ is a smooth solution of the following third order ordinary differential equation:
\begin{equation}\label{SC}\begin{cases}
\Big(\tilde\tau+sign(U_s''(a))z^2\Big)^2\frac{d}{dz}W+i\frac{d^3}{dz^3}\Big(\big(\tilde\tau+sign(U_s''(a))z^2\big)W\Big)=0,\\
\lim\limits_{z\rightarrow-\infty}W~=~0,\quad\lim\limits_{z\rightarrow+\infty}W~=~1.
\end{cases}\end{equation}

The approximate solution $(u_\ep^{app},v_\ep^{app})(t,x,y)$ given in \eqref{app_1} can be  used to prove
 the instability of the problem \eqref{sim_lin} because the expression \eqref{app_1} combining with
 the property of the
parameter $\tau: \Im\tau<0$  implies a growing mode $e^{-\frac{\Im\tau}{\sqrt{\ep}}}$ 
for the approximation  $(u_\ep^{app},v_\ep^{app})(t,x,y)$ when $\ep\ll1$. 
However, plugging the formula \eqref{for_v} of $v_\ep(y)$ into \eqref{app_1}   yields 
\begin{equation}\label{decay_u}
u_\ep^{app}(t,x,y)~=~e^{i\ep^{-1}(x+w_\ep t)}v_\ep'(y),\quad v_\ep'(y)~=~H(y-a)U_s'(y)+\ep^{\frac{1}{4}}V'(\frac{y-a}{\ep^{\frac{1}{4}}}).
\end{equation}
Then, it implies that the approximation $u_\ep^{app}(t,x,y)$ has the same decay rate as $U_s'(y)$ when $\yinf$. In particular, $u_\ep^{app}\notin \h^m_\alpha$ initially for any $\alpha>0$ if $U_s'(y)$ does not decay exponentially as $\yinf$.

 Therefore, to study the case of shear flow with general decay, 
the above  approximation $u_\ep^{app}$ in \eqref{decay_u} 
will be inappropriate since the operator we consider now is 
$$T(t,s): ~\h^{m_1}_\alpha~\mapsto~\h_0^{m_2},\quad \forall \alpha>0,~\mbox{for some}~m_1,m_2>0.$$
For this, we need to modify the construction of approximate solution \eqref{app_1} with \eqref{for_v} to problem \eqref{sim_lin}, in order that at least the initial tangential data of the approximation has an exponential decay rate as $\yinf.$ So, we 
will look for a new approximate solution of \eqref{sim_lin} in the following form:
\begin{equation}\label{app_2}
(\tilde u_\ep^{app},\tilde v_\ep^{app})(t,x,y)~=~e^{i\ep^{-1}(x+\tilde w_\ep t)}\Big(iv_{\ep,1}'(y)+itv_{\ep,2}'(y),\ep^{-1}v_{\ep,1}(y)+\ep^{-1} tv_{\ep,2}(y)\Big).
\end{equation}
In the above expression, we expect that, on one hand,
\begin{equation}\label{ass_app}
\tilde w_\ep~=~w_\ep,\quad v_{\ep,2}(y)~=~v_\ep(y),
\end{equation}
where $w_\ep$ and $v_\ep(y)$ are given in Proposition \ref{prop_GD},
thus the instability of \eqref{sim_lin} preserves through the eigenvalue perturbation $\tau$ as mentioned above; on the other hand,
\begin{equation}\label{bd_tv}
v_{\ep,1}(0)~=~v_{\ep,1}'(0)~=~0,\qquad \lim\limits_{y\rightarrow+\infty}v_{\ep,1}'(y)=0,~\mbox{exponentially},
\end{equation}
so that the initial data of $\tilde u^{app}_\ep(t,x,y)$ given by \eqref{app_2} has an exponential decay rate as $\yinf$.

The motivation of the construction  in the form
of \eqref{app_2} comes from the expression of solutions to the linearized inviscid Prandtl equation around a shear flow
$\big(U(y),0\big)$. That is, the system
\begin{equation}\label{pr_invis}\begin{cases}
\pd_t u+U(y)\pd_x u+U'(y)v=0,\\
\pd_x u+\pd_y v=0,\\
v|_{y=0}=0,\quad u|_{t=0} =u_0(x,y)
\end{cases}\end{equation}
has the solution
\begin{equation}\label{express_invis}\begin{split}
&u(t,x,y)=u_0\big(x-tU(y),y\big)+tU'(y)\int_0^y u_{0x}\big(x-tU(z),z\big)dz,\\
&v(t,x,y)=-\int_0^y\Big\{u_{0x}\big(x-tU(z),z\big)+t\big[U(y)-U(z)\big] u_{0xx}\big(x-tU(z),z\big)\Big\}dz,
\end{split}\end{equation}
 see Proposition 5.1 in \cite{H-H}. From the above expression \eqref{express_invis}, we know that when $t>0$, the decay rate of tangential velocity of
the solution to the  problem \eqref{pr_invis} is not faster than the one of background shear flow $U'(y)$, even though the initial data $u_0(x,y)$ can decay very rapidly as $\yinf$.

Now, it remains to find a suitable $v_{\ep,1}(y)$ for the new approximation \eqref{app_2}, such that the error term
\[
\tilde r^{app}_\ep~:=~\pd_t \tilde u^{app}_\ep+\cl_s\tilde u^{app}_\ep
\]
still satisfies the relation \eqref{est_err1}. By virtue of \eqref{ass_app}, a direct computation yields that $\tilde r^{app}_\ep(t,x,y)=e^{i\ep^{-1}(x+ w_\ep t)}\tilde R^{app}_\ep(y)$ and
\begin{equation}\label{for_err}\begin{split}
\tilde R^{app}_\ep(y)~&=~-\ep^{-1}[w_\ep+U_s(y)]v_{\ep,1}'(y)+\ep^{-1}U_s'(y)v_{\ep,1}(y)+iv_{\ep,1}^{(3)}(y)+iv_\ep'(y)+t R^{app}_\ep(y)\\
~:&=~\bar R_\ep^{app}(y)+t R_\ep^{app}(y)
\end{split}\end{equation}
with $R^{app}_\ep(y)$ given by \eqref{for_R}.
Note that
\[
-\ep^{-1}[w_\ep+U_s(y)]v_{\ep,1}'(y)+iv_{\ep,1}^{(3)}(y)\in W_\alpha^{0,\infty}(\R^+)
\]
provided that $v_{\ep,1}'(y)\in W_\alpha^{0,\infty}(\R^+)$ for some $\alpha>0.$ 
Thus, to ensure $\bar R^{app}_\ep(y)\in W_\alpha^{0,\infty}(\R^+)$, we only need
\[
\ep^{-1}U_s'(y)v_{\ep,1}(y)+iv_\ep'(y)\in W_\alpha^{0,\infty}(\R^+),
\]
which implies that by combining with \eqref{decay_u}, 
\begin{equation}\label{decay_u1}
v_{\ep,1}(y)~\rightarrow~-i\ep\quad\mbox{exponentially,}\quad\mbox{as}~\yinf.
\end{equation}
Obviously, for any function $f(y),~y\in\R^+$:
\begin{equation}\label{def_f}
	f(y)\in C_c^\infty(\R^+),\qquad \int_0^{+\infty}f(y)dy\neq0,
\end{equation}
the function
\begin{equation}\label{for_v1}
v_{\ep,1}(y)~:=~-i\ep\frac{\int_0^yf(z)dz}{\int_0^{+\infty}f(y)dy}
\end{equation}
meets the requirements \eqref{bd_tv} and \eqref{decay_u1}. Then, plugging the above expression \eqref{for_v1} into \eqref{app_2},  we obtain  the new approximate solution to \eqref{sim_lin}:
\begin{equation}\label{app_3}
(\tilde u_\ep^{app},\tilde v_\ep^{app})(t,x,y)~=~e^{i\ep^{-1}(x+ w_\ep t)}\Big(\frac{ \ep f(y)}{\int_0^{+\infty}f(y)dy} +itv_{\ep}'(y),-\frac{i\int_0^yf(z)dz}{\int_0^{+\infty}f(y)dy}+\ep^{-1} tv_{\ep,2}(y)\Big),
\end{equation}
where the functions  $v_\ep(y)$ and $f(y)$ are given by \eqref{for_v} and \eqref{def_f} respectively.

\subsection{Construction of approximate solutions}

Following the construction of approximate solutions to the simplified problem \eqref{sim_lin} given in
the previous subsection, and also by the arguments  used
in \cite{GV-D}, we are going to construct the approximate solutions to the original linearized problem \eqref{linear_pr}. Since the approximate solutions to \eqref{sim_lin}  given in \eqref{app_3}  are obtained with the background state being frozen at the initial data $u_s|_{t=0}=U_s(y)$, to construct the approximate solutions of the original problem \eqref{linear_pr} with background state being shear flow in the time interval $0<t<t_0$, we need  some modification as in \cite{GV-D}.

Let $u_s(t,y)$ satisfy the assumptions of Theorem \ref{thm_lin},
and $a>0$ be a non-degenerate critical point of $U_s(y)$. Without loss of generality, we assume that $U_s''(a)<0$, then the differential equation
\begin{equation}\label{def_alpha}\begin{cases}
\pd_t\pd_yu_s\big(t,a(t)\big)+\pd_y^2u_s\big(t,a(t)\big)a'(t)=0,\\
a(0)~=~a
\end{cases}\end{equation}
defines a non-degenerate critical point $a(t)$ of $u_s(t,\cdot)$ when $0<t<t_0$ for some small $t_0>0$. Moreover, we have $\pd_y^2u_s\big(t,a(t)\big)<0$ for all $t\in[0,t_0)$ with $t_0$ small enough.
As in \cite{GV-D}, we take $\tau, W(z)$ given by \eqref{SC} (we drop the tilde of $\tilde{\tau}$ for brevity), and set
\begin{equation}\begin{split}\label{def_phi}
V(z)~&:=~\big(\tau-z^2\big) W( z)-1_{\R^+}\big(\tau-z^2\big).
\end{split}\end{equation}
For $0<\ep\ll1$, introduce
\begin{equation}\label{app_lam}
w_\ep(t)~:=~-u_s\big(t,a(t)\big)+\ep^{\frac{1}{2}}\Big|\frac{\pd_y^2u_s\big(t,a(t)\big)}{2}\Big|^{\frac{1}{2}}~\tau,
\end{equation}
and the ``regular'' part of the tangential velocity field
\begin{equation}\label{app_reg}
v_\ep^{reg}(t,y)=H\big(y-a(t)\big)\Big[u_s(t,y)-u_s\big(t,a(t)\big)+\ep^{\frac{1}{2}}\Big|\frac{\pd_y^2u_s\big(t,a(t)\big)}{2}\Big|^{\frac{1}{2}} \tau\Big],
\end{equation}
as well as the ``shear layer'' part
\begin{equation}\label{app_sl}\begin{split}
v_\ep^{sl}(t,y)~&:=~\ep^{\frac{1}{2}}\varphi\big(y-a(t)\big)\Big|\frac{\pd_y^2u_s\big(t,a(t)\big)}{2}\Big|^{\frac{1}{2}} ~V\Big(\Big|\frac{\pd_y^2u_s\big(t,a(t)\big)}{2}\Big|^{\frac{1}{4}}\cdot\frac{y-a(t)}{\ep^{\frac{1}{4}}}\Big).
\end{split}\end{equation}
Here, $\varphi$ is a smooth truncation function near 0, and $V(z)$ is given in \eqref{def_phi}. Also, for any function $f(y)$ satisfying \eqref{def_f},
let
\begin{equation}\label{for-tv}
\tilde v_{\ep}(y)~:=~\frac{\int_0^yf(z)dz}{\int_0^{+\infty}f(y)dy}.
\end{equation}

Next, according to the discussion in the above subsection, the approximate solution of the problem \eqref{linear_pr} can be defined as follows:
\begin{equation}\label{ap_lin}
(u_\ep,v_\ep)(t,x,y)=e^{i\ep^{-1}x}\big(U_\ep,V_\ep\big)(t,y)
\end{equation}
with
\begin{equation}\begin{split}\label{expre_uv}
U_\ep(t,y)~&=~e^{i\ep^{-1}\int_0^tw_\ep(s)ds}\Big[\ep\tilde v'_\ep(y)+it\pd_y\big(v_\ep^{reg}(t,y)+v_\ep^{sl}(t,y)\big)\Big],\\
V_\ep(t,y)~&=~e^{i\ep^{-1}\int_0^tw_\ep(s)ds}\Big[-i\tilde v_\ep(y)+\ep^{-1}t\big( v_\ep^{reg}(t,y)+v_\ep^{sl}(t,y)\big)\Big].
\end{split}\end{equation}
For the function $(u_\ep,v_\ep)(t,x,y)$ in \eqref{ap_lin} to be $2\pi-$periodic in $x$, we take $\ep=\frac{1}{n}$ with  $n\in\N$. It is straightforward to check that, 
\[(u_\ep,v_\ep)|_{y=0}=0,\qquad\limz u_\ep=0,\]
 and the divergence free condition holds. Also, $u_\ep(t,x,y)=e^{i\ep^{-1}x}U_\ep(t,y)$ is analytic in the tangential variable $x$ and $W^{2,\infty}$ in $y$. Moreover, there are positive constants $C_0$ and $\sigma_0$, independent of $\ep$, such that
\begin{equation}\label{bound_app}
\|U_\ep(t,\cdot)\|_{W_0^{2,\infty}}
~\leq~C_0e^{\frac{\sigma_0t}{\sqrt{\ep}}},\quad t\in[0,t_0),
\end{equation}
 in particular,
\begin{equation}\label{bound_ini} \|U_\ep(0,\cdot)\|_{W_\alpha^{2,\infty}}
~\leq~C_0 \ep,\qquad\forall\alpha\geq0.
\end{equation}

Plugging the relation \eqref{ap_lin} into the original linearized Prandtl equations \eqref{linear_pr}, it follows that
\begin{equation}\label{lin_ap}\begin{cases}
\pd_t u_\ep+u_s\pd_xu_\ep+v_\ep \pd_yu_s-\pd_y^2u_\ep=r_\ep,\quad&\\
\pd_xu_\ep+\pd_yv_\ep=0,\quad&{\rm in}\quad \Omega,\\
(u_\ep,v_\ep)|_{y=0}=0.
\end{cases}\end{equation}
The remainder term $r_\ep$ can be represented by $r_\ep(t,x,y)=e^{i\ep^{-1}x}R_\ep(t,y)$
and
\begin{equation}\label{r}\begin{split}
R_\ep(t,y)~:=~\bar R_\ep(t,y)+t\tilde R_\ep(t,y),
\end{split}\end{equation}
where
\begin{equation}\label{br}\begin{split}
\bar R_\ep(t,y)=e^{i\ep^{-1}\int_0^tw_\ep(s)ds}\Big\{&i\big[w_\ep(t)+u_s(t,y)\big]\tilde v'_\ep(y)-i\pd_yu_s(t,y)\tilde v_\ep(y)-\ep \tilde v_\ep^{(3)}(y)\\
&+i\pd_y\big(v_\ep^{reg}(t,y)+v_\ep^{sl}(t,y)\big)\Big\},
\end{split}\end{equation}
and
\begin{equation}\label{tr}\begin{split}
\tilde R_\ep(t,y)=e^{i\ep^{-1}\int_0^tw_\ep(s)ds}\Big\{&-\ep^{-1}\Big[u_s(t,y)
-u_s(t,a(t))-\pd_y^2u_s(t,a(t))\frac{(y-a(t))^2}{2}\Big]\pd_yv_\ep^{sl}(t,y)\\
&+\ep^{-1}\Big[\pd_yu_s(t,y)-\pd_y^2u_{s}(t,a(t))(y-a(t))\Big]v_\ep^{sl}(t,y)\\
&+i\pd_t\pd_yv_\ep^{sl}(t,y)+O(\ep^\infty)\Big\}.
\end{split}\end{equation}
The term $O(\ep^\infty)$ in \eqref{tr} represents the part of remainder
with exponential decay in $y$ that comes from the fact that $V(z)$ decays exponentially and the derivatives of $\varphi(\cdot-a(t))$ vanish 
outside a neighborhood of  $a(t)$. Combining the formulation \eqref{tr} of $\tilde R_\ep(t,y)$ and the exponential decay of $v_\ep^{sl}(t,y)$ yields
\begin{equation}\label{bd_tr}
\|\tilde R_\ep(t,\cdot)\|_{W_\alpha^{0,\infty}}\leq C_1e^{\frac{\sigma_0t}{\sqrt{\ep}}}, \quad\forall\alpha\geq0
\end{equation}
with the constant $\sigma_0>0$ given in \eqref{bound_app}. On the other hand,
 from  \eqref{app_reg}-\eqref{for-tv} we have
\[
-i\pd_yu_s(t,y)\tilde v_\ep(y)+i\pd_y\big(v_\ep^{reg}(t,y)+v_\ep^{sl}(t,y)\big)~\equiv~0,\quad \mbox{for large}~ y>0,
\]
and then,
\[
-i\pd_yu_s(t,y)\tilde v_\ep(y)+i\pd_y\big(v_\ep^{reg}(t,y)+v_\ep^{sl}(t,y)\big)\in W_\alpha^{2,\infty},\quad\forall\alpha\geq0,
\]
which implies that the estimate \eqref{bd_tr} also holds for the term $\bar R_\ep(t,y)$.
Thus,  with the same $\sigma_0$ given in \eqref{bound_app}, the term
$R_\ep(t,y)$ satisfies
\begin{equation}\label{bound_r}
\|R_\ep(t,\cdot)\|_{W_\alpha^{0,\infty}}\leq C_1e^{\frac{\sigma_0t}{\sqrt{\ep}}}, \quad\forall\alpha\geq0,
\end{equation}
where the constant $C_1>0$ is independent of $\ep$.

\subsection{Proof of the main Theorem}

Based on the approximate solutions construted in the above subsection, we can 
apply the approach  in \cite{GV-D} to prove  Theorem \ref{thm_lin}. We
now sketch the proof as follows.

The proof is based on the verification of \eqref{est_in} for the tangential differential operator by contradiction.
Suppose that \eqref{est_in} does not hold, that is, for all $\sigma>0$, there exists $\delta>0,~\alpha_0,m\geq0$ and $\mu\in[0,\frac{1}{2})$ such that
\begin{equation}\label{pr_x}
\sup\limits_{0\leq s< t\leq \delta}\|e^{-\sigma(t-s)\sqrt{|\pd_x|}}T(t,x)\|_{\cl(\h_{\alpha_0}^m,\h_0^{m-\mu})}<+\infty.
\end{equation}
Introduce the operator
\begin{equation*}
T_\ep(t,s):~W_{\alpha_0}^{0,\infty}(\R^+)\mapsto W_{0}^{0,\infty}(\R^+)
\end{equation*}
as
\begin{equation}\label{def_t}
 T_\ep(t,s)U_0~:=~e^{-i\ep^{-1}x}T(t,s)\Big(e^{i\ep^{-1}x}U_0\Big)
\end{equation}
with $T(t,s)$ being defined in \eqref{def_T}.
From \eqref{pr_x}, we have
\begin{equation}\label{est_tep}
\|T_\ep(t,s)\|_{\cl(W_{\alpha_0}^{0,\infty},W_0^{0,\infty})}~\leq~C_2\ep^{-\mu}e^{\frac{\sigma(t-s)}{\sqrt{\ep}}},\quad \forall ~0\leq s< t\leq\delta
\end{equation}
for some constant $C_2>0$ independent of $\ep$.

Next, denote by
\[L_\ep~:=~e^{-i\ep^{-1}x}~L~e^{i\ep^{-1}x},\]
where $L$ is the linearized Prandtl operator around the shear flow $\big(u^s(t,y),0\big)$.
Let $U(t,y)$ be a solution to the problem
\[
\pd_tU+L_\ep U~=~0,\qquad
U|_{t=0}~=~U_\ep(0,y),
\]
where $U_\ep(t,y)$ is given in \eqref{expre_uv}.
Thus, we have 
\[U(t,y)~=~T_\ep(t,0)U_\ep(0,y),\]
and by using \eqref{bound_ini} and \eqref{est_tep} it follows that 
\begin{equation}\label{up_bound}
\|U(t,\cdot)\|_{W_0^{0,\infty}}\leq C_2\ep^{-\mu}e^{\frac{\sigma t}{\sqrt{\ep}}}\|U_\ep(0,\cdot)\|_{W_{\alpha_0}^{0,\infty}}
\leq C_3\ep^{1-\mu}e^{\frac{\sigma t}{\sqrt{\ep}}},\quad \forall t\in(0,\delta]
\end{equation}
 for some constant $C_3>0$ independent of $\ep$.

On the other hand, we know that the difference $\tilde U:=U-U_\ep$ can be obtained by the Duhamel principle:
\begin{equation}\label{eq_dif}
\tilde U(t,\cdot)~=~\int_0^tT_\ep(t,s) R_\ep(s,\cdot)ds,\quad \forall~t\leq\delta.
\end{equation}
From \eqref{bound_r}, \eqref{est_tep} and \eqref{eq_dif},
and choosing $\sigma<\sigma_0$, we have
\begin{equation}\label{est_dif}
\|\tilde U(t,\cdot)\|_{W^{0,\infty}_0}\leq C_1C_2\ep^{-\mu}
\int_0^te^{\frac{\sigma(t-s)}{\sqrt{\ep}}}e^{\frac{\sigma_0s}{\sqrt{\ep}}}ds
\leq C_4\ep^{\frac{1}{2}-\mu}e^{\frac{\sigma_0t}{\sqrt{\ep}}},
\end{equation}
where the constant $C_4>0$ is independent of $\ep$. Then, by combining  \eqref{est_dif} with the expression  of $U_\ep(t,y)$ in \eqref{expre_uv},
we obtain that for  $t\in(0,\delta]$ and sufficiently small $\ep$,
\begin{equation}\label{low_bound}\begin{split}
\|U(t,\cdot)\|_{W_0^{0,\infty}}&\geq\|U_\ep(t,\cdot)\|_{W_0^{0,\infty}}
-\|\tilde U(t,\cdot)\|_{W_0^{0,\infty}}\\
&\geq e^{\frac{\sigma_0t}{\sqrt{\ep}}}(C_5 t-C_6\ep)-C_4\ep^{\frac{1}{2}-\mu}e^{\frac{\sigma_0t}{\sqrt{\ep}}}\\
&\geq C_5te^{\frac{\sigma_0t}{\sqrt{\ep}}}-2C_4\ep^{\frac{1}{2}-\mu}e^{\frac{\sigma_0t}{\sqrt{\ep}}}.
\end{split}\end{equation}
As $\sigma<\sigma_0,$ comparing \eqref{up_bound} with \eqref{low_bound},
the contradiction arises when $t\gg\frac{\mu|\ln \ep|}{\sigma_0-\sigma}\ep^{\frac{1}{2}-\mu}$
with sufficiently small $\ep$. Thus, the proof of Theorem \ref{thm_lin} is completed.

\section{Further discussions}
In this section, we  point out that the above results can be extended to the three space dimensions under some condition on the background shear flow given in \cite[Theorem 2.3]{LWY2}. More precisely, consider the three dimensional Prandtl equations in the domain $\{(t,x,y,z):~t>0,(x,y)\in\T^2,z\in\R^+\}:$
\begin{equation}\label{3dpd}
\begin{cases}
\pd_t u+(u\pd_x+v\pd_y+w\pd_z)u-\pd_z^2u=0,\\
\pd_t v+(u\pd_x+v\pd_y+w\pd_z)v-\pd_z^2v=0,\\
\pd_xu+\pd_yv+\pd_zw=0,\\
(u,v,w)|_{z=0}=0,\qquad \lim\limits_{z\rightarrow+\infty} (u,v)=(U_0,V_0)
\end{cases}\end{equation}
with positive constants $U_0$ and $V_0$.
Let $(u_s,v_s)(t,z)$ be a smooth solution of the heat equations:
\begin{equation}\label{shear}\begin{cases}
\pd_t u_s-\pd_z^2u_s=0,\qquad
\pd_t v_s-\pd_z^2v_s=0,\\
(u_s,v_s)|_{z=0}=0,\qquad \lim\limits_{z\rightarrow+\infty}  (u_s,v_s)=(U_0,V_0).
\end{cases}\end{equation}
It is
straightforward to verify that the shear velocity profile $(u_s,v_s,0)(t,z)$ solves the problem \eqref{3dpd}.
Then, we study the linearized problem of \eqref{3dpd} around the shear flow $(u_s,v_s,0)(t,z)$:
\begin{equation}\label{lin}\begin{cases}
\pd_t u+(u_s\pd_x+v_s\pd_y)u+w\pd_zu_s-\pd_z^2u=0,\\
\pd_t v+(u_s\pd_x+v_s\pd_y)v+w\pd_zv_s-\pd_z^2v=0,\\
\pd_xu+\pd_yv+\pd_zw=0,\\
(u,v,w)|_{z=0}=0,\qquad\lim\limits_{z\rightarrow+\infty} (u,v)=0.
\end{cases}\end{equation}

Denote by $T(t,s)$ the linearized solution operator of problem \eqref{lin}, i.e.,
\begin{equation}\label{def_T}
T(t,s)\big((u_0,v_0)\big)~:=~(u,v)(t,\cdot),
\end{equation}
where $(u,v)$ is the solution of \eqref{lin} with $(u,v)|_{t=s}=(u_0,v_0)$.  The result on the linear instability of the three-dimensional Prandtl equations is:

\begin{prop}\label{prop_3d}
	Let $(u_s,v_s)(t,z)$ slove \eqref{shear} with
	$$(u_s-U_0,v_s-V_0)\in C^0\big(\R^+; W_0^{4,\infty}(\R_z^+)\big)\bigcap C^1\big(\R^+; W_0^{2,\infty}(\R_z^+)\big),$$
	and assume that the initial data $(U_s,V_s)(z)\triangleq(u_s,v_s)(0,z)$ satisfies that
	\begin{equation}\label{ass}
	\exists~ z_0>0,~s.t.~
	V_s'(z_0)U_s''(z_0)\neq U_s'(z_0)V_s''(z_0).
	\end{equation}
	Then, 
	 there exists $\sigma>0$ such that for any $\delta>0$,
	\begin{equation}
	\sup\limits_{0\leq s< t\leq\delta}\big\|e^{-\sigma(t-s)\sqrt{|\pd_\cT|}}T(t,s)\big\|_{\cl(H_\alpha^m,H_0^{m-\mu})}
	~=~+\infty,\quad \forall m,\alpha\geq0,~\mu\in[0,\frac{1}{4}),
	\end{equation}
	where the operator $\pd_\cT$ represents the tangential derivative $\pd_x$ or $\pd_y$, and the weighted Sobolev spaces $H_\alpha^m$ are given by
	\begin{equation*}
	H_\alpha^m~:=~H^m\big(\T^2_{x,y};W_\alpha^{0,\infty}(\R^+_z)\big),\qquad\forall m,\alpha\geq0.
	\end{equation*}
	Moreover, 
	\begin{equation}
	\sup\limits_{0\leq s< t\leq\delta}\big\|T(t,s)\big\|_{\cl(H_\alpha^m,H_0^{0})}
	~=~+\infty,\quad \forall m,\alpha\geq0.
	\end{equation}

\end{prop}
This proposition can be proved by combining the above arguments with the analysis in  \cite{LWY2}, hence,  we omit it for brevity.

Finally, the nonlinear instability in both
2D and 3D cases can also be discussed for
the case when the background shear flow has general decay by using 
the above linear
instability results and the arguments from \cite{GV-N,guo} and \cite{LWY2}.

\vspace{.15in}

{\bf Acknowledgement:}
The  research  was supported by the General Research Fund of Hong Kong,
CityU No. 103713.


\end{document}